\documentclass[reqno, 11pt, a4paper, oneside]{amsart}
\usepackage{amsmath}
\usepackage{amsfonts}
\usepackage{amssymb}
\usepackage[ansinew]{inputenc}
\usepackage{graphicx}
\usepackage{amsbsy}
\usepackage{fancyhdr}
\usepackage{yfonts}
\usepackage{hyperref}
\usepackage{mathabx}
\usepackage{enumerate}

\numberwithin{equation}{section}
\newtheorem{teo}{Theorem}[section]
\newtheorem{prop}[teo]{Proposition}
\newtheorem{lema}[teo]{Lemma}

\theoremstyle{definition}
\newtheorem{defi}[teo]{Definition}

\theoremstyle{remark}

\sffamily
\title{Local uniformity through larger scales}
\author{Miguel N. Walsh}
\address{Departamento de Matemática e IMAS-CONICET, Facultad de Ciencias Exactas y Naturales, Universidad de Buenos Aires, 1428 Buenos Aires, Argentina}
\email{mwalsh@dm.uba.ar}

\begin{document}

\def\F{\mathbb{F}}
\def\Fqn{\mathbb{F}_q^n}
\def\Fq{\mathbb{F}_q}
\def\Fp{\mathbb{F}_p}
\def\Di{\mathbb{D}}
\def\E{\mathbb{E}}
\def\Z{\mathbb{Z}}
\def\Q{\mathbb{Q}}
\def\C{\mathbb{C}}
\def\R{\mathbb{R}}
\def\N{\mathbb{N}}
\def\P{\mathbb{P}}
\def\Pc{\mathcal{P}}
\def\T{\mathbb{T}}
\def\modp{\, (\text{mod }p)}
\def\modN{\, (\text{mod }N)}
\def\modq{\, (\text{mod }q)}
\def\modone{\, (\text{mod }1)}
\def\Zn{\mathbb{Z}/N \mathbb{Z}}
\def\Zp{\mathbb{Z}/p \mathbb{Z}}
\def\Zan{a^{-n}\mathbb{Z}/ \mathbb{Z}}
\def\Zal{a^{-l} \Z / \Z}
\def\Pr{\text{Pr}}
\def\leftsize{\left| \left\{}
\def\rightsize{\right\} \right|}

\maketitle

\begin{abstract}
By associating frequencies to larger scales, we provide a simpler way to derive local uniformity of multiplicative functions on average from the results of Matomäki-Radziwi\l\l.
\end{abstract}

\section{Introduction}

In this note we present a simpler derivation of the following estimate.

\begin{teo}
\label{1}
Given $0 < \delta < 1$ and $\eta > 0$, there exists some $C>0$ such that, for every complex-valued multiplicative function $g$ with $|g| \le 1$ satisfying
\begin{equation}
\label{A}
 \int_X^{2X} \sup_{\alpha} \left| \sum_{x \le n \le x+H} g(n) e(\alpha n) \right| dx \ge \eta H X,
 \end{equation}
with $H=X^{\delta}$, we have $\mathbb{D} (g; CX^2/H^2,C) \le C$.
\end{teo}

Here, as usual, $\mathbb{D}(g;T,Q)$ stands for the 'pretentious' distance introduced by Granville and Soundararajan \cite{GS}:
$$ \mathbb{D}(g;T,Q) = \inf \left( \sum_{p \le T} \frac{1 - \text{Re}(g(p)p^{it}\chi(p))}{p} \right)^{1/2} ,$$
with the infimum taken over all $|t| \le T$ and all Dirichlet characters of modulus at most $Q$.

Theorem \ref{1} was originally obtained by Matomäki, Radziwi\l\l \text{ }and Tao \cite{MRT}. Their original result had an $H^{\varepsilon}$-loss but this was later removed by Matomäki, Radziwi\l\l, Tao, Teräväinen and Ziegler \cite{ALL} in the more general setting of nilsequences. The result constitutes progress towards what is termed the local uniformity conjecture which in the linear case states that Theorem \ref{1} should hold assuming only that $H \rightarrow_X \infty$. Besides being a natural question, it was shown by Tao \cite{T} that its restriction to the Liouville function becomes equivalent to both Chowla's and Sarnak's conjectures after performing a logarithmic averaging. A number of other consequences also follow from this conjecture (see \cite{ALL}).

It should be straightforward to adapt the approach in this paper to the case of polynomial phases. Similarly, the structure of the argument, which reduces matters to an application of Vinogradov's lemma, would seem well-suited for an appeal to the equidistribution theory of nilsequences (see \cite{GT}) and this may facilitate the task of improving the known range of validity of the full local uniformity conjecture and its consequences. However, in the present article we shall focus solely in deriving Theorem \ref{1}, which showcases neatly the main idea.

We proceed through the same setting as in \cite{MRT}. For any given $x$ there are essentially a bounded number of choices of $\alpha$ for which the integrand in (\ref{A}) can be big. For simplicity, let us assume in this discussion there is at most one such choice. By Elliot's inequality and the multiplicativity of $g$, this means that we can find a lot of pairs $x,x' \in [X,2X]$ and primes $p,p'$ in some range $[P,2P]$ with $x/p$ close to $x'/p'$ and with their corresponding frequencies $\alpha, \alpha'$ satisfying that $\alpha p$ is close to $\alpha' p'$. From this information, our goal is to show that $g(n)$ must 'pretend' to be $n^{2\pi iT} \chi(n)$ for some $T=O(X^2 /H^2)$ and a certain Dirichlet character $\chi$ of bounded modulus.

To accomplish this we look at intervals of the form $[y,y+PH]$ with $y \sim PX$ and we show that the compatibility between the frequencies that we just described allows us to associate to many of these intervals a frequency $\alpha_y$ for which we can find a lot of primes $p \in [P,2P]$ such that if $y/p$ is close to $x$ and $\alpha$ was the frequency associated to the interval $[x,x+H]$, then $p \alpha_y$ is close to $\alpha$. This in turn implies that these new elements $\alpha_y$ and intervals $[y,y+PH]$ satisfy the same kind of relations that we had over the original intervals and phases. As a consequence, we can apply an iterative procedure that ultimately leads us to a point $z \sim X^2/H$ and a corresponding frequency $\alpha_z$ such that for many of the original intervals $[x,x+H]$, we can find products of primes $p_1 \cdots p_k$ with $\frac{z}{p_1 \cdots p_k}$ close to $x$ and $p_1 \cdots p_k \alpha_z$ close to the frequency $\alpha$ corresponding to this interval. This establishes a 'global' relationship between the original frequencies and through an application of Vinogradov's lemma, implies that they must essentially correspond to a modulation of the first order Taylor approximation of $n^{2 \pi i T}$ for an adequate $T$. The claim then follows from the results in \cite{MR,MRT2}.

The strategy just described is relatively easy to implement and allows us to avoid a lot of the machinery used in \cite{MRT,ALL}, which includes the need to build up a large modulus for the relationships between the frequencies, the use of primes at different scales, graph-theoretic and probabilistic arguments, a mixing lemma based on the Vinogradov-Korobov bound and the use of approximate ergodicity. This may have further advantages due to known limitations of some of these previous tools (see the end of \cite[Section 1.1]{MRT2}).

The idea of studying the 'residual' behaviour of $g(n)$ on larger scales would seem to be useful for problems that involve averaged estimates. One could think of the frequencies we construct as phases $g(n)$ would correlate with shall we extend it multiplicatively to these new scales. An inspection of the argument going back to the proof of Lemma \ref{conc} below also provides some intuition on the kind of 'conspiracies' that may occur to prevent the validity of the conjectural bounds: for an element $y \sim PX$ like we described before, the points $x$ that are close to $y/p$ for some $p \in [P,2P]$ should group into a bounded number of subsets such that the frequency associated to a given $x$ satisfies the compatibility conditions with essentially all elements in the same subset but practically none of those in other subsets. Hopefully both the idea of working through larger scales and the rigidity just described could prove relevant in future work on this and related problems.

\subsection{Notation} We will write $X \lesssim Y$ or $X=O(Y)$ to mean that there is some absolute constant $C$ with $|X| \le C Y$ and $X \sim Y$ if both $X \lesssim Y$ and $Y \lesssim X$ hold. If the implicit constants depend on some additional parameters, we shall use a subscript to indicate this. For a finite set $S$ we write $|S|$ for its cardinality. We shall write ${\bf 1}_{P(n)}$ to refer to the indicator function of those elements $n$ satisfying a certain property $P(n)$. We abbreviate $e(x):=e^{2 \pi i x}$ and write $\| \cdot \|$ for the distance to $0$ in $\R / \Z$. Finally, given an interval $J \subseteq \R$, we write $J/p$ for the set of integers $m$ with $mp \in J$.

\section{Contagiousness of phase relations}

Our starting point is the following simple observation.

\begin{lema}
\label{start}
Let $\epsilon > 0$. Let $\alpha_1,\alpha_2 \in \R / \Z$ and let $p_1,p_2$ be primes with $\| p_1 \alpha_2 - p_2 \alpha_1 \| < \epsilon$. Then, there exists $\alpha \in \R / \Z$ with $\| p_i \alpha - \alpha_i \| < \frac{\epsilon}{2p_j}$ if $i \neq j$.
\end{lema}

\begin{proof}
For $i=1,2$, let $\alpha^{(i)} \in \R / \Z$ be such that $p_i \alpha^{(i)} = \alpha_i$. Adding integer multiples of $1/p_i$ to $\alpha^{(i)}$ we may assume that $\| \alpha^{(1)} - \alpha^{(2)} \| \le \frac{1}{2 p_1 p_2}$. On the other hand, we have by hypothesis that $\| p_1 p_2 (\alpha^{(1)} - \alpha^{(2)} ) \| < \epsilon$. Combining both estimates we see that it must in fact be $\| \alpha^{(1)} - \alpha^{(2)} \| < \frac{\epsilon}{p_1 p_2}$. Taking $\alpha$ to be the closest middle point of $\alpha^{(1)}$ and $\alpha^{(2)}$ in $\R / \Z$, we obtain the result.
\end{proof}

These kind of relations can easily be seen to be contagious in the following sense.

\begin{lema}
\label{tg}
Let $\alpha_1,\alpha_2,\alpha_3 \in \R / \Z$. Let $P > 1$ and let $0 \le \epsilon < \frac{c}{P}$ for some sufficiently small $c \gtrsim 1$. Let $P \le p_1,p_2,p_3 \le 2P$ be distinct primes with $\| p_i \alpha_j - p_j \alpha_i \| < \epsilon$ for every choice of $i,j \in \left\{ 1,2,3 \right\}$. Let $\alpha$ be the element of $\R / \Z$ associated to $\alpha_1,\alpha_2,p_1,p_2$ by Lemma \ref{start}. Then $\| p_3 \alpha - \alpha_3 \| =O( \frac{\epsilon}{P})$.
\end{lema}

\begin{proof}
Arguing as in the proof of Lemma \ref{start}, for each pair $i, j \in \left\{ 1,2,3 \right\}$ we can find elements $\alpha_{\{ i,j \}}^{(i)},\alpha_{\{ i,j \}}^{(j)} \in \R / \Z$ with $p_k \alpha_{\{ i,j \}}^{(k)} = \alpha_k$ for $k \in \left\{ i,j \right\}$, $ \left\| \alpha_{\{ i,j \}}^{(i)}-\alpha_{\{ i,j \}}^{(j)} \right\|<  \frac{\epsilon}{p_i p_j} < \frac{\epsilon}{P^2}$ and such that $\alpha$ is the closest middle point of $\alpha_{\{ 1,2 \}}^{(1)}$ and $\alpha_{\{ 1,2 \}}^{(2)}$. But since $\alpha_{\{ i,j \}}^{(i)}-\alpha_{\{ i,k \}}^{(i)}=\frac{l_i}{p_i}$ for some integer $0 \le l_i < p_i$, we see from several applications of the triangle inequality that 
$$ \frac{\epsilon}{P^2} >  \left\| \alpha_{\{ 1,2 \}}^{(1)} -\alpha_{\{ 1,2 \}}^{(2)} \right\| =  \| \frac{l_1}{p_1} + \frac{l_2}{p_2}+\frac{l_3}{p_3} \| + O(\frac{\epsilon}{P^2}  ).$$
It follows that $l_1=l_2=l_3=0$, for otherwise the first term would have magnitude at least $(2P)^{-3}$ and this would contradict the size of $\epsilon$. In particular, $\alpha_{\{ 1,2 \}}^{(1)} =\alpha_{\{ 1,3 \}}^{(1)}$. Since $\| \alpha- \alpha_{\{ 1,2 \}}^{(1)} \| = O(\frac{\epsilon}{P^2})$, we deduce from the triangle inequality that $\| \alpha - \alpha_{\{ 1,3 \}}^{(3)} \| = O( \frac{\epsilon}{P^2})$ and the result immediately follows.
\end{proof}

Similarly, we have the following variant in which we have two pairs with both elements of each pair relating with every element of the other pair.

\begin{lema}
\label{2}
Let $P > 1$ and let $0 \le \epsilon < \frac{c}{P^2}$ for some sufficiently small $c \gtrsim 1$. Let $P \le p_1,p_2,q_1,q_2 \le 2P$ be distinct primes and let $\alpha_{p_1},\alpha_{p_2},\alpha_{q_1},\alpha_{q_2} \in \R / \Z$. Suppose that $\| p_i \alpha_{q_j} - q_j \alpha_{p_i} \| < \epsilon$ for every choice of $1 \le i,j \le 2$. Then we also have $\| p_1 \alpha_{p_2} - p_2 \alpha_{p_1} \| \lesssim \epsilon$.
\end{lema}

\begin{proof}
As in the proof of Lemma \ref{start}, given a choice of  $p_i,q_j$ we can find elements $\beta_j^{(i)},s_j^{(i)}  \in \R / \Z$ with $\| s_j^{(i)} \| <  \frac{\epsilon}{P^2}$ such that $p_j \beta_j^{(i)}= \alpha_{p_j}$ and $q_i (\beta_j^{(i)}+s_j^{(i)}) = \alpha_{q_i}$. The latter implies that $ \beta_1^{(i)}-\beta_2^{(i)} = \frac{l^{(i)}}{q_i} + O(\frac{\epsilon}{P^2})$ for some integer $0 \le l^{(i)} < q_i$. Similarly, we have that  $ \beta_j^{(1)}-\beta_j^{(2)} = \frac{l_j}{p_j}$ for some integer $0 \le l_j < p_j$. Therefore
$$ \frac{l^{(1)}}{q_1} + O(\frac{\epsilon}{P^2})=  \beta_1^{(1)}-\beta_2^{(1)} = \frac{l_1}{p_1}  - \frac{l_2}{p_2} + \beta_1^{(2)} - \beta_2^{(2)} = \frac{l_1}{p_1}  - \frac{l_2}{p_2} + \frac{l^{(2)}}{q_2} + O(\frac{\epsilon}{P^2}) ,$$
implying that
$$ \| \frac{l_1}{p_1}  - \frac{l_2}{p_2} + \frac{l^{(2)}}{q_2}-\frac{l^{(1)}}{q_1} \| = O(\epsilon/P^2).$$
If the left side were nonzero, it would have magnitude at least $(2P)^{-4}$, contradicting the size of $\epsilon$. This implies that each term must vanish. We have thus found elements $\beta_1, \beta_2 \in \R / \Z$ with $p_j \beta_j = \alpha_{p_j}$ and $\| \beta_1 - \beta_2 \| =O( \epsilon/P^2)$. Multiplying $\beta_1 - \beta_2$ by $p_1 p_2$, we conclude the proof.
\end{proof}

The preceding lemmas can be combined to yield the following estimate, which is what we shall need later.

\begin{lema}
\label{conc}
Let $P > 1$ and let $0 \le \epsilon < \frac{c}{P^2}$ for some sufficiently small $c \gtrsim 1$. Let $S \subseteq [P,2P]$ be a set of primes. For each $p \in S$, let $\alpha_p$ be a corresponding element of $\R / \Z$. If for at least $\gtrsim |S|^2$ of the pairs $p_1,p_2 \in S$ we have $\| p_1 \alpha_{p_2} - p_2 \alpha_{p_1} \| < \epsilon$, then there exists an element $\alpha \in \R / \Z$ such that $\| p \alpha - \alpha_p \| \lesssim \frac{ \epsilon}{P}$ for $\gtrsim |S|$ primes $p \in S$. 
\end{lema}

\begin{proof}
We may assume $|S| > C$ for a sufficiently large constant $C$. For each $p \in S$ write $S_p$ for the set of $p' \neq p \in S$ satisfying $\| p \alpha_{p'} - p' \alpha_p \| < \epsilon$. Since there are $\gtrsim |S|^2$ pairs in $S$ satisfying this relation, we may find distinct primes $p_1,p_2$ with  $|S_{p_1} \cap S_{p_2}| \gtrsim |S|$. By Lemma \ref{2} applied with some pair $q_1 \neq q_2 \in S_{p_1} \cap S_{p_2}$ it must be $\| p_1 \alpha_{p_2} - p_2 \alpha_{p_1} \| \lesssim \epsilon$ and by Lemma \ref{start} we can hence find some $\alpha \in \R / \Z$ such that $\| p_i \alpha - \alpha_{p_i} \| \lesssim \frac{\epsilon }{P} $ for $i=1,2$. The result now follows from applying Lemma \ref{tg} to $p_1,p_2,p_3$ for each $p_3 \in S_{p_1} \cap S_{p_2}$.
\end{proof}

\section{Setting up the recursion}

The goal of this section is mostly to recreate the setting in \cite{MRT}, as described in the introduction. The exception is given by Proposition \ref{it}, which is concerned with the structure of the iterative procedure we shall apply.

\begin{defi}
Given $H > 1$, $c > 0$ and an interval $I$, we say $\mathcal{J} \subseteq I \times \R / \Z$ is a $(c,H)$-configuration in $I$ if $|\mathcal{J}| \ge c|I|/H$ and the first coordinates are $H$-separated points in $I$ (i.e. $|x-y| \ge H$ if $x \neq y$).
\end{defi}

As in \cite{MRT}, it is convenient for us to discretise the problem. The following lemma to this effect is immediate.

\begin{lema}
\label{el}
Suppose Theorem \ref{1} fails and let the notation be as in its statement. Then there exist a $(c_0,H)$-configuration $\mathcal{J}_1$ in $[X,2X]$, for some $c_0 \gtrsim \eta$, such that for every $(x,\alpha_x) \in \mathcal{J}_1$ it is
\begin{equation}
\label{la}
 \left| \sum_{x < n \le x+H} g(n) e(\alpha_x n) \right| \ge \eta H/2.
 \end{equation}
\end{lema}

To build the relations between the frequencies we shall use the following form of Elliot's inequality (see for example \cite[Proposition 2.5]{MRT}).

\begin{lema}
\label{elliot}
Let $\tau >0$. Let $f:I \rightarrow \C$ with $|f(n)| \le 1$ for every $n$ in some interval $I$. Then
$$ \frac{1}{|I|} \sum_{n \in I} f(n) - \frac{p}{|I|} \sum_{n \in I/p} f(n) \lesssim \tau,$$
for all primes $p \le |I|$ outside of an exceptional set of primes $\mathcal{P}$ satisfying $\sum_{p \in \mathcal{P}} \frac{1}{p} \lesssim \tau^{-2}$.
\end{lema}

Applying this in our setting, we obtain the following lemma.

\begin{lema}
\label{q}
Let the notation and assumptions be as in Lemma \ref{el} and let $\varepsilon > 0$ be sufficiently small with respect to $\eta$ and $\delta$. Then, there exists some integer $H^{\varepsilon^2} \lesssim P \lesssim H^{\varepsilon}$ such that, for at least $\gtrsim |\mathcal{J}_1|$ elements $(x,\alpha_x) \in \mathcal{J}_1$, there are at least $\gtrsim \frac{P}{\log P}$ primes $p \in [P,2P]$ with
\begin{equation}
\label{Jp}
\left| \frac{p}{H}\sum_{m \in [x+1,x+H]/p} g(pm) e(\alpha_x pm) \right| \ge \eta/4.
 \end{equation}
\end{lema}

\begin{proof}
Given $(x,\alpha_x) \in \mathcal{J}_1$ and a prime $p$, we write $p \dashv (x,\alpha_x)$ to mean that (\ref{Jp}) is satisfied. Write $P_0 = \lfloor H^{\varepsilon^2} \rfloor$ and recursively define $P_i = 2^i P_0$. Let $k$ be the smallest integer with $P^{k} > H^{\varepsilon}$. Notice that if $\varepsilon > 0$ is chosen sufficiently small with respect to $\eta$, we can guarantee that 
$$\log \log P_k - \log \log P_0 \ge \log \log H^{\varepsilon} - \log \log H^{\varepsilon^2} \ge C \eta^{-2},$$   
for some sufficiently large constant $C$. It thus follows from Lemma \ref{el} and Mertens' theorem that
\begin{equation}
\begin{aligned}
 \sum_{(x,\alpha_x) \in \mathcal{J}_1} \sum_{i=0}^k \sum_{p \in [P_i,2P_i)} \frac{{\bf 1}_{p \dashv (x,\alpha_x)}}{p} &\gtrsim |\mathcal{J}_1| (\sum_{P_0 \le p \le 2 P_k} \frac{1}{p}- O(\eta^{-2})) \\
 &\gtrsim  |\mathcal{J}_1| \frac{\log \log P_k - \log \log P_0}{2}. 
 \end{aligned}
 \end{equation}
The result is now an easy consequence of the prime number theorem and the pigeonhole principle.
\end{proof}

The next lemma clarifies the claim made in the introduction that the integrand in (\ref{A}) can be big for essentially only a bounded number of $\alpha$.

\begin{lema}
\label{LS}
Let $\tau > 0$. Let $g: I \rightarrow \C$ be defined on an interval $I$ with $|g(n)| \le 1$ for every $n \in I$. Then, there exists a set $S \subseteq \R / \Z$ with $|S|=O_{\tau}(1)$ such that for every $\alpha \in \R / \Z$ with
\begin{equation}
\label{ab}
\left| \frac{1}{|I|} \sum_{n \in I} g(n) e(\alpha n) \right| \ge \tau,
\end{equation}
we have $\| \alpha - \beta \| \lesssim_{\tau} |I|^{-1}$ for some $\beta \in S$.
\end{lema}

\begin{proof}
Let $\alpha_1,\ldots,\alpha_K$ be distinct elements of $\R / \Z$ satisfying (\ref{ab}). By a standard application of the Cauchy-Schwarz inequality, we can find complex coefficients $c_1,\ldots,c_K$ with $\sum_{i=1}^K |c_i|^2 = 1$ such that
$$ \sum_{1 \le i,j \le K} c_i \overline{c_j}\sum_{n \in I} e(\alpha_i n-\alpha_j n) \gtrsim K \tau^2 |I|.$$
The result follows noticing that if $K$ is sufficiently large with respect to $\tau$ as to make the non-diagonal contribution dominate, we can find a pair $i \neq j$ with
$$\tau^2 |I| \lesssim \sum_{n \in I} e((\alpha_i - \alpha_j )n) \lesssim \frac{1}{\| \alpha_i - \alpha_j\|} .$$
\end{proof}

Using this input and the multiplicativity of $g$, we can now start to build the relations between the frequencies.

\begin{lema}
\label{base}
Let the notation and assumptions be as in Lemma \ref{q} and assume $X$ is sufficiently large with respect to $\eta$ and $\delta$. Then, there exists some $c \gtrsim_{\eta} 1$ and a $(c,H/P)$-configuration $\mathcal{J}_0$ in $[X/(2P),2X/P]$ such that, for $\gtrsim_{\eta} X/H$ elements $(x,\alpha_x)$ of $\mathcal{J}_1$ we can find $\gtrsim_{\eta} \frac{P}{\log P}$ pairs $(p,(y,\alpha_y)) \in [P,2P] \times \mathcal{J}_0$ with $p$ prime, such that $|py -x| < 2H \text{ and } \| p \alpha_x - \alpha_y \| \lesssim_{\eta} P/H$.
\end{lema}

\begin{proof}
By Lemma \ref{q} we know that for $\gtrsim |\mathcal{J}_1|$ elements $(x,\alpha_x) \in \mathcal{J}_1$, there are at least $\gtrsim \frac{P}{\log P}$ primes $p \in [P,2P]$ satisfying (\ref{Jp}). Notice that $g(pm)=g(p)g(m)$ for all except perhaps the $O(H/P^2)$ terms of the sum (\ref{Jp}) with $p|m$. Therefore, in each case where (\ref{Jp}) is satisfied we can take a subinterval $I \subseteq [x+1,x+H]/p$ of length $\lfloor H/2P \rfloor$ with
$$ \left| \frac{2P}{H}\sum_{m \in I} g(m) e(\alpha_x p m)) \right| \ge \eta/8,$$
say. Furthermore, any other interval $I'$ of the same length with initial point $y$ and with symmetric difference $| I \Delta I' | \le \frac{\eta H}{32 P}$ will satisfy the same estimate with $\eta/16$ in place of $\eta/8$ and by a simple application of the triangle inequality we see that it will also obey $|py-x|<2H$. Considering all pairs $(y,p \alpha_x)$ generated in this way for some $P \le p \le 2P$ and some $(x,\alpha_x) \in \mathcal{J}_1$ and applying the pigeonhole principle, we can find a set $ \mathcal{J}_0^*$ of  $H/P $-separated points in $[X/(2P),2X/P]$, with $|\mathcal{J}_0^*| \gtrsim_{\eta} X/H$, each of whose elements relates in the above manner with $\gtrsim_{\eta} \frac{P}{\log P}$ distinct $(x,\alpha_x) \in \mathcal{J}_1$. Applying Lemma \ref{LS} we see that for each of the $\gtrsim_{\eta} \frac{P}{\log P}$ elements $(x,\alpha_x) \in \mathcal{J}_1$ with which a given $y \in \mathcal{J}_0^*$ relates in this way, we have $\| p \alpha_x - \beta_i \| \lesssim_{\eta} P/H$ for some $1 \le i \le r$, where $\beta_1,\ldots,\beta_r$ is a fixed set of $r=O_{\eta}(1)$ elements of $\R / \Z$ that we can associate to $y$. In particular, there will be one element in this set that will satisfy this for $\gtrsim_{\eta} \frac{P}{\log P}$ elements of $\mathcal{J}_1$. Writing $\alpha_y$ for this element and double-counting, the set of pairs $(y,\alpha_y) \in \mathcal{J}_0^* \times \R / \Z$ gives us the desired configuration $\mathcal{J}_0$.
\end{proof}

The next proposition shows how to use the kind of relations we have obtained to build similar relations at a larger scale.

\begin{prop}
\label{it}
Let $A,P > 1$, $c_1,c_2 \sim 1$, $k \in \N$ and $H > CP^2$ for some sufficiently large $C \sim 1$. Let $\mathcal{P}$ be the set of primes in $[P,2P]$. Let $\mathcal{J}_k$ be a $(c_1,HP^k)$-configuration in an interval $[Y, AY]$ and let $\mathcal{J}_{k+1}$ be a $(c_2,HP^{k+1})$-configuration in $[PY,2 P A Y]$. Assume that for each $y \in \mathcal{J}_{k+1}$ there are at least $\gtrsim \frac{P}{\log P}$ pairs $(p,(x,\alpha_x)) \in \mathcal{P} \times \mathcal{J}_k$ such that $ |px-y| \lesssim HP^{k+1}$ and $\| p \alpha_{y} - \alpha_{x} \| \lesssim (HP^{k} )^{-1}$. Then, there exists $c_3 \gtrsim 1$ and a $(c_3,HP^{k+2})$-configuration $\mathcal{J}_{k+2}$ in $[P^2 Y, (2P)^2 AY]$ such that, for every $(z,\alpha_z) \in \mathcal{J}_{k+2}$, there are at least $\gtrsim \frac{P}{\log P}$ pairs $(p,(y,\alpha_y)) \in \mathcal{P} \times \mathcal{J}_{k+1}$ with $ |py-z| \lesssim HP^{k+2}$ and $\| p \alpha_{z} - \alpha_{y} \| \lesssim (HP^{k+1})^{-1}$.
\end{prop}

\begin{proof}
Let $(x,\alpha_x) \in \mathcal{J}_k$ and suppose $p_1,p_2 \in \mathcal{P}$ satisfy $|p_i x-y_i| \lesssim HP^{k+1} \text{ and } \| p_i \alpha_{y_i} - \alpha_{x} \| \lesssim (HP^k)^{-1} $, where $(y_i,\alpha_{y_i})$ are elements of $\mathcal{J}_{k+1}$. Multiplying by $p_i$ and using the triangle inequality, it follows that $|p_2 y_1 - p_1 y_2| \lesssim HP^{k+2}$ and  $\| p_1 \alpha_{y_1} -p_2 \alpha_{y_2} \| \lesssim (HP^k)^{-1}$. Combining our assumptions with an application of the Cauchy-Schwarz inequality, we see that the number of quintuples $(x,y_1,y_2,p_1,p_2)$ of the above form is $\gtrsim \frac{X}{H} \frac{P^{2}}{(\log P)^{2}}$. Notice also that for each such quintuple we have that $p_i y_j \in [P^2 Y, (2P)^2 AY]$ and that given $p \in \mathcal{P}$, if $|y-y'| \ge HP^{k+1}$, then $|py - py'| \ge HP^{k+2}$. It thus follows easily from the pigeonhole principle that we can locate a set $\mathcal{J}_{k+2}^{\ast}$ of $HP^{k+2}$-separated points in $[P^2 Y, (2P)^2 AY]$ of cardinality $\gtrsim X/H$ such that, if for each $z \in \mathcal{J}_{k+2}^{\ast}$ we write $r (z)$ for the number of quadruples $(p_1,p_2,(y_1,\alpha_{y_1}),(y_2,\alpha_{y_2})) \in (\mathcal{P})^2 \times \mathcal{J}_{k+1}^2 $ with $|z - p_i y_i| \lesssim HP^{k+2}$ and $\| p_2 \alpha_{y_1}-p_1\alpha_{y_2} \| \lesssim (HP^k)^{-1}$ (notice how the roles of $p_1$ and $p_2$ have been inverted with respect to the reasoning we did before), we have
$$ \sum_{z \in \mathcal{J}_{k+2}^{\ast}} r(z) \gtrsim  \frac{X}{H} \frac{P^{2}}{(\log P)^{2}}.$$
We may discard from $\mathcal{J}_{k+2}^*$ those $z$ with $r(z) \le c \frac{P^{2}}{(\log P)^{2}}$ without affecting the estimate as long as $c$ is chosen sufficiently small. For each of the $\gtrsim X/H$ remaining $z$ we obtain from Lemma \ref{conc} an element $\alpha_z \in \R / \Z$ such that the set $\mathcal{J}_{k+2}$ of pairs $(z,\alpha_z)$ gives us the desired configuration.
\end{proof}

We will need two further lemmas. The first one is the following estimate on products of primes \cite[Lemma 2.6]{MRT}.

\begin{lema}
\label{products}
Let $k,Q \in \N$ and let $P,N \ge 3$ with $P^{k-1} \gtrsim N$. Then, the number of tuples of primes $(p_1,\ldots,p_k,p_1',\ldots,p_k') \in [P,2P]^{2k}$ with
$$ \left| \prod_{j=1}^k p_i - \prod_{j=1}^k p_i' \right| \lesssim P^k/N,$$
and with both products congruent mod $Q$, is at most $\lesssim_k \frac{P^{2k}}{N(\log P)^{2k}} (1/\varphi(Q) + 1/\log N)$.
\end{lema}

Finally, we will also need the following standard estimate \cite[Lemma 1.1.14]{T0}.

\begin{lema}[Vinogradov's lemma]
\label{v}
Let $\alpha \in \R / \Z$, let $0 < \epsilon < 1/100$, $100 \epsilon < \delta < 1$ and $\delta N > 100$. Suppose that $\| n \alpha \| < \epsilon$ for at least $\delta N$ integers $n \in [-N,N]$. Then, there exist integers $a,q \lesssim \delta^{-1}$ with $\| \alpha - a/q \| \lesssim \frac{\epsilon}{\delta N}$.
\end{lema}

\section{The proof of Theorem \ref{1}}

Let the notation and assumptions be as in Lemma \ref{base}. From now on all implicit constants are allowed to depend on $\eta$ and $\delta$. We may assume $X>C$ for a sufficiently large constant $C$. Let $\tilde{k}_0$ be the smallest positive integer with $\frac{P^{\tilde{k}_0}}{(\log P)^{\tilde{k}_0+1}} > X/H$ and let $\tilde{k}=\tilde{k}_0+1$. In particular, $\tilde{k} \lesssim 1$. Notice that given $p,q \in [P,2P]$, if $|py-z| \lesssim HP^{k+2}$ and $\| p \alpha_{z} - \alpha_{y} \| \lesssim (HP^{k+1})^{-1}$ and we also have $|qx-y| \lesssim HP^{k+1}$ and $\| q \alpha_{y} - \alpha_{x} \| \lesssim (HP^{k})^{-1}$, then an application of the triangle inequality shows that $|pqx-z| \lesssim HP^{k+2}$ and $\| pq \alpha_{z} - \alpha_{x} \| \lesssim (HP^{k})^{-1}$. As a consequence, we see that iterating $\tilde{k}$ times Proposition \ref{it} using Lemma \ref{base} as a starting point, we end up with a $(c,HP^{\tilde{k}})$-configuration $\mathcal{J}_{\tilde{k}}$ in $[XP^{\tilde{k}}, AXP^{\tilde{k}}]$, for some $c,A \sim 1$, such that for each $(z,\alpha_z) \in \mathcal{J}_{\tilde{k}}$ there exist $\gtrsim \frac{P^{\tilde{k}}}{(\log P)^{\tilde{k}}}$ elements $(x,\alpha_x) \in \mathcal{J}_1$, counting repetitions, for which we can find primes $p_1, \ldots, p_{\tilde{k}} \in [P,2P]$ with $|p_1 \cdots p_{\tilde{k}} x - z| \lesssim HP^{\tilde{k}}$ and $\| p_1 \cdots p_{\tilde{k}} \alpha_{z} - \alpha_{x} \| \lesssim H^{-1}$. Notice that by passing to a subset of $\mathcal{J}_{\tilde{k}}$ of size $\gtrsim X/H$ if necessary and adjusting the implicit constants accordingly, we may assume that a pair $((z,\alpha_z),p_1 \cdots p_{\tilde{k}})$ can be associated in the above way to at most one element of $\mathcal{J}_1$.

It follows that given $(z,\alpha_z) \in \mathcal{J}_{\tilde{k}}$ we can find an element $(x,\alpha_x) \in \mathcal{J}_1$ such that the above relation holds for at least $\gtrsim \frac{H}{X} \frac{P^{\tilde{k}}}{(\log P)^{\tilde{k}}}$ products $p_1 \ldots p_{\tilde{k}}$. Labelling these products $q_1,\ldots,q_m$ we see from the triangle inequality that for every $1 \le j \le m$ we have $(q_j-q_1) x \lesssim HP^{\tilde{k}}$, so $(q_j-q_1) \lesssim \frac{HP^{\tilde{k}}}{X}$, and also $\| (q_{j}- q_1)\alpha_{z} \| \lesssim H^{-1}$. An application of Lemma \ref{v} then tells us that $\alpha_z = \frac{a}{Q} + O ( \frac{ X (\log P)^{\tilde{k}}}{H^2 P^{\tilde{k}}} )$ for some coprime integers $a,Q \lesssim (\log P)^{\tilde{k}}$. This immediately implies that for each element $(x,\alpha_x) \in \mathcal{J}_1$ satisfying the relations of the above paragraph with $(z,\alpha_z)$, the corresponding products $p_1 \cdots p_{\tilde{k}}$ will have to be congruent mod $Q$. But we already know that we have $\gtrsim \frac{P^{\tilde{k}}}{(\log P)^{\tilde{k}}}$ such products across all elements of $\mathcal{J}_1$. An application of Cauchy-Schwarz would then contradict Lemma \ref{products} with $N=X/H$ unless $Q=O(1)$ and at least $\gtrsim X/H$ distinct elements $\mathcal{J}_1^{(1)} \subseteq \mathcal{J}_1$ relate with $(z,\alpha_z)$ via at least $\gtrsim \frac{H}{X} \frac{P^{\tilde{k}}}{(\log P)^{\tilde{k}}}$ products. We now let $z_1,\ldots,z_r$ be $r \sim (\log P)^{\tilde{k}}$ points that are $\gtrsim \frac{H P^{\tilde{k}}}{(\log P)^k}$-separated such that given any $(x,\alpha_x) \in \mathcal{J}_1$ and product $p_1 \cdots p_{\tilde{k}}$ that relate with $(z,\alpha_z)$ like before, we have $|p_1 \cdots p_{\tilde{k}} x - z_i| \lesssim \frac{HP^{\tilde{k}}}{(\log P)^{\tilde{k}}}$ for some $1 \le i \le r$. In particular, given any two such products of primes $q,q'$ relating like this with the same $z_i$, we have $|q-q'| \lesssim \frac{HP^{\tilde{k}}}{X (\log P)^{\tilde{k}}}$. Let $\epsilon \gtrsim 1$ be sufficiently small and let $\mathcal{J}_1^{(2)} \subseteq \mathcal{J}_1^{(1)}$ consist of those elements for which there are less than $\epsilon (\log P)^k$ choices of $1 \le i \le r$ for which we can find some product $p_1 \cdots p_{\tilde{k}}$ that relates with $z_i$ in this manner. Appealing again to Cauchy-Schwarz and Lemma \ref{products}, now with $N= \frac{X (\log P)^{\tilde{k}}}{H}$, we see that it must be $\mathcal{J}_1^{(1)} \setminus \mathcal{J}_1^{(2)} \gtrsim X/H$. This in turn allows us to find, by the pigeonhole principle, some $1 \le i \le r$ such that for $\gtrsim X/H$ elements $(x,\alpha_x) \in \mathcal{J}_1$ we have $|p_1 \cdots p_{\tilde{k}} x - z_i| \lesssim \frac{HP^{\tilde{k}}}{(\log P)^{\tilde{k}}}$ and $\| p_1 \cdots p_{\tilde{k}} \alpha_{z} - \alpha_{x} \| \lesssim H^{-1}$ for some product $p_1 \cdots p_{\tilde{k}}$ of primes in $[P,2P]$. It also implies that we can find two such products of primes $q,q'$ with $|q- q'| \sim \frac{HP^{\tilde{k}}}{X}$ and $\|(q-q')\alpha_z \| \lesssim H^{-1}$, which allows us to improve our rational approximation of $\alpha_z$ to $\alpha_z= \frac{a}{Q} + O(\frac{X}{H^2 P^{\tilde{k}}})$.

Since $z_i \sim X P^{\tilde{k}}$, we can thus write $ \alpha_z = \frac{a}{Q} + \frac{T}{z_i}$, for some real number $T \lesssim X^2/H^2$. Furthermore, for each of the $\gtrsim X/H$ elements $(x,\alpha_x) \in \mathcal{J}_1$ that relate with $z_i$ via a product $p_1 \cdots p_{\tilde{k}}$, we have by the triangle inequality that
\begin{equation}
\label{aprox}
 \| \alpha_x - (\frac{p_1 \cdots p_{\tilde{k}} a}{Q} + \frac{T}{x}) \| \le \| \frac{T}{x} - \frac{p_1 \cdots p_{\tilde{k}} T}{z_i} \|  + O(H^{-1}) \lesssim H^{-1}.
 \end{equation}
We can now finish the proof of Theorem \ref{1} as in the end of \cite[Section 5]{MRT}. If $\epsilon \gtrsim 1$ is sufficiently small, then for any of the $\gtrsim X/H$ elements $(x,\alpha_x) \in \mathcal{J}_1$ obeying (\ref{la}) and (\ref{aprox}) we have that
$$ \sum_{y < n < y+H_{\ast}} g(n) e( \frac{p_1 \cdots p_{\tilde{k}} an}{Q} + \frac{T(n-y)}{y}) \gtrsim H_{\ast},$$
if $|y - x| \le H_{\ast}:=\epsilon H$. Considering the first order Taylor approximation of $n^{2 \pi iT}$, performing a Fourier decomposition of $e( \frac{p_1 \cdots p_{\tilde{k}} an}{Q})$ into Dirichlet characters and applying the pigeonhole principle, we can thus find a Dirichlet character of modulus $O(1)$ with
$$ \sum_{y < n < y+H_{\ast}} g(n) n^{2 \pi iT} {\bf 1}_{r | n} \chi(n/r) \gtrsim H_{\ast},$$
for $\gtrsim X$ intervals $[y,y+H_{\ast}] \subseteq [X,2X]$ and some fixed $r \sim 1$. But this can easily be seen to contradict \cite[Theorem A.1]{MRT2} unless $\mathbb{D}(f; O(T),O(1)) \lesssim 1$, concluding the proof of Theorem \ref{1}.

\end{document}